\documentclass[11pt,english]{amsart}
\usepackage{hyperref}
\usepackage[left=1in,top=1in,right=1in,bottom=1in]{geometry}
\usepackage{wasysym, mathpazo, mathtools, amssymb, tikz, mathrsfs}
\pdfoutput=1
\usetikzlibrary{matrix,arrows}

\input xy 
\xyoption{all}

\newtheorem{theorem}{Theorem}[section]

\newtheorem{conjecture}[theorem]{Conjecture}
\newtheorem{proposition}[theorem]{Proposition}
\newtheorem{corollary}[theorem]{Corollary}

\theoremstyle{definition}
\newtheorem{definition}[theorem]{Definition}
\newtheorem{example}[theorem]{Example}

\theoremstyle{remark}
\newtheorem{remark}[theorem]{Remark}

\newcommand{\Aff}{{\mathbb A}}
\newcommand{\C}{{\mathbb C}}
\newcommand{\PP}{{\mathbb P}}
\newcommand{\Q}{{\mathbb Q}}

\newcommand{\Z}{{\mathbb Z}}
\newcommand{\Qbar}{{\overline{\Q}}}

\def\bbar#1{\setbox0=\hbox{$#1$}\dimen0=.2\ht0 \kern\dimen0 
\overline{\kern-\dimen0 #1}}
\newcommand{\Kbar}{{\bbar{K}}}  

\newcommand{\calC}{{\mathcal C}}
\newcommand{\calO}{{\mathcal O}}

\DeclareMathOperator{\supp}{supp}
\DeclareMathOperator{\Pic}{Pic}

\newcommand{\tcal}[1]{{\tilde{\mathcal{#1}}}}
\newcommand{\Oint}[1]{{\calO_{#1}}} 

\begin{document}
\title{Invitation to Integral and Rational points on curves and surfaces.} 
\author{Pranabesh Das}
\address{Stat-Math Unit, Indian Statistical Institute, 7 S.J.S. Sansanwal Marg, New Delhi 110 016, India}
\email{pranabesh.math@gmail.com} 
\author{Amos Turchet}
\address{Universit\`{a} degli Studi di Udine, Dipartimento di Matematica e Informatica, via delle Scienze 208, 33100 Udine, Italy}
\curraddr{Department of Mathematical Sciences, Chalmers University of Technology and the University of Gothenburg, SE-412 96 Gothenburg, Sweden}

\email{tamos@chalmers.se}
\date{}

\maketitle
\section{Introduction}\label{sec:intro}
Solving polynomial equations with integral coefficients is one of the oldest mathematical problems and has been studied, in its easiest formulation, since the ancient times. These equations are known at the present time under the name of \emph{Diophantine Equations}. This title is owed to the Greek mathematician Diophantus (AD 200), who wrote the book \emph{Arithmetica} which contains several problems regarding solutions of polynomial equations. The subject grew in importance through all mathematical history and, in the seventeenth century, was put in the spotlight by Pierre de Fermat and his famous claim that the equation
\[
x^n+y^n=z^n \qquad n\geq3
\]
has no solution with $x,y,z$ non-zero rationals. Fermat found this interesting problem while reading Diophantus' book, but unfortunately his copy had a too narrow margin to contain the proof (``\textit{Hanc marginis exiguitas non caperet}'')\footnote{``Whether Fermat knew a proof or not has been the subject of many speculations. The truth seems to be obvious (...). Fermat quickly became aware of the incompleteness of the "proof" of 1637. Of course, there was no reason for a public retraction of his privately made conjecture.'' \cite{FermatMinkowski}}. The proof was given by Andrew Wiles (building on works of many others, among whom we cite Hellegouarch, Frey, Ribet and Taylor) more than three hundred years after Fermat's note. In the journey for the search of a solution, new subjects were born, like Algebraic Number Theory, and numerous topics, results and techniques were developed to attack this particular equation, e.g. ideals, UFD, Class Number, Dirichlet series, L-functions.

Among these results, a special mention should be given to Falting's proof of Mordell's Conjecture which shows, in particular, that for all $n\geq 3$ the number of solutions to Fermat's Equations, if there exist any, is finite. This breakthrough was a milestone in the field of \emph{Diophantine Geometry}, a subject whose name comes from the seminal first edition of Lang's book \cite{Lang1983}. The idea behind this subject is to study the properties of the set of solutions to polynomial equations over an algebraic closed field containing $\Q$, usually $\C$ or $\Qbar$, instead of $\Q$. Such a set of solutions has a geometric structure, which turns it into a so-called \emph{algebraic variety}, and its geometric properties can give insights on the number and the distribution of solutions to the original equation.

This expository article focus on the Diophantine problems related to equations in three or four variables, which correspond geometrically to the arithmetic study of curves and surfaces. After recalling some basic facts and techniques about Diophantine Equations, to which is devoted Section 2, we describe the distribution of integral and rational points on algebraic curves in Section 3, via the theorems of Siegel and Faltings. In the last section we address some of the recent developments in (complex) dimension 2 focusing on two deep conjectures stated, respectively, by Bombieri and (independently) Lang, and by Vojta.

No particular prerequisites in Algebraic Geometry or Number Theory are required for the first three sections: we only assume some elementary facts about Arithmetic and easy geometric constructions accessible to any undergraduate student. For the last section, instead, we need some techniques coming from Algebraic Geometry at a basic level, such as \cite{Hartshorne} or \cite{Shafarevich}. We are not claiming any originality in the results presented here, except for some computations and examples. Our main goal is to give an informal introduction to this fascinating subject, hoping to show to the interested reader the beauty of these results and, rephrasing Diophantus, ``\textit{the nature and power subsisting in numbers}''. 

\section{Some Diophantine Equations}\label{sec:DioEq}
In his well-known speech \cite{Hilbert1900} given at the 1900 International Congress of Mathematicians in Paris, David Hilbert posed the following problem, which has become known as \textit{Hilbert 10th problem}:

\begin{quote}
Given a Diophantine equation with any number of unknown quantities\\
and with rational integral numerical coefficients to devise a process\\
according to which it can be determined by a finite number of\\
operations whether the equation is solvable in rational integers.
\end{quote}

In 1970 Matiyasevish gave a definite negative answer, dropping the hopes to obtain an algorithm which determines a priori whether a polynomial equation has solutions or not (see \cite{Matiassevitch} for a detailed exposition). However, this does not mean that classes of Diophantine equations cannot be solved, even explicitly. In this section we are going to see some elementary methods for dealing with particular equations that will be important in the subsequent sections. The methods used here are elementary in nature but still provides the answer to the questions we are interested in: "Does there exists a solution to this equations in either rationals or integrals? If any, how many solutions are there?". We begin by giving an explicit example which resembles Fermat's equation:

\begin{example}\label{ex:ratsols1}
Consider the following quadratic equation in three variables
\begin{equation}\label{eq1}
x^2+y^2=z^2.
\end{equation}
Suppose we want to find all the rational solutions to equation (\ref{eq1}). Immediately, we can see that putting all the variables to zero gives a trivial solution. However, this is not interesting because (\ref{eq1}) is homogeneous and it always has a zero solution. The solutions we are really interested in are the ones which represent a point of the quadric defined by (\ref{eq1}) in $\PP^2$. Indeed, to each (non-zero) rational solution $(x,y,z)$ in $\Q^3$, by multiplying each coordinate by a non zero constant we could obtain infinitely many solutions. Instead, we should be looking at solutions which are \emph{projective points}, i.e. triples $\left[ x : y : z \right]$ defined up to multiplication by non-zero scalars. In particular $x,y,$ and $z$ cannot be all zero, removing the trivial solution found above.

From this observation it follows that there are no projective rational solutions with $z=0$, since this would imply that both $x$ and $y$ were 0. Then, we can assume that $z\neq 0$ and divide both sides by $z^2$; in this way we get an equation of the form
\[
X^2+Y^2=1,
\]
which represents a circle in the affine plane $\Aff^2$, where $X$ and $Y$ are $x/z$ and $y/z$, respectively. The question of finding rational solutions to equation (\ref{eq1}) now becomes a matter of finding rational solutions to this new equation.

Note that $(-1,0)$ is an integral, and therefore rational, solution; we claim that from such solution we can find infinitely many others. Let us look at all lines with rational slopes passing through (-1,0), say $y=mx+m$ with $m$ rational. When we plug the line equation into the circle equation we get a quadratic equation
\[
(m^2+1)x^2+2m^2x+(m^2-1)=0,
\] 
whose solutions are
\[
x_m =  \dfrac{(-2m^2\pm{2})}{2(m^2+1)}.
\]
Now, this line cuts the circle in a second point different from $(-1,0)$, namely
\[
P_m = \Big(\frac{1-m^2}{1+m^2},\frac{2m}{1+m^2}\Big),
\]
which is rational and, for all $m$, all these points are different from each other. This proves that the equation $x^2+y^2=z^2$ has infinitely many rational solutions. Even more, here we proved that all rational (projective) solutions are in a one to one correspondence with the lines passing through $(-1,0)$. See Example \ref{ex:dioph_ex} for a similar problem solved with a geometric argument.
\end{example}

Above we discussed what might be seen as an ad-hoc method to solve one particular equation. In this next example we are going to see how the strategy can be actually applied to many others quadratic equations in three variables:

\begin{example}\label{ex:ratsols2}
Consider the following Diophantine equation:
\begin{equation}\label{eq2}
2x^2+z^2=3y^2.
\end{equation}
By the same argument we sketched before we divide both sides by $z^2$ and this time we get the equation of an hyperbola:
\[
2 X^2 +1 = 3 Y^2.
\]
Now, let us consider the problem of finding all rational solutions to this new equation. Again we can easily find one solution, i.e. $(1,1)$. Then, we look at the lines $y=mx-(m-1)$ with $m$ rational, i.e. lines with rational slope passing through $(1,1)$. Then, by the same method of Example \ref{ex:ratsols1}, we get a quadratic equation. This time the discriminant is $4 (2 - 3 m)^2$ which is a square for all rational $m$. Therefore, there exists infinitely many $m$'s which gives a rational solution to the equation in one variable. This implies that there are infinitely many rational solutions both to the hyperbola equation and to equation (\ref{eq2}). In particular, every time we can find a rational solution to the affine equation obtained from the original homogeneous equation, the method of using lines with rational slopes passing through the particular solution gives infinitely many rational solutions of the equation we are interested in. The only requirement for this method to work is that the equation one gets by plugging the line equation into the conic one has infinitely many rational solutions. This translate the problem into a number theoretic one, i.e. into discussing whether the discriminant is a square for some values of the slopes. However, since the lines have rational slopes and intersects the conic in one rational point one can prove that the other point of intersection is again a point with rational coordinates.
\end{example}

Previous example can be generalized in the following:

\begin{theorem}
Let us consider a homogeneous polynomial equation in three variables of degree two with integral or rational coefficients. If there exists a non-zero rational solution then there exist infinitely many rational solutions.
\end{theorem}

The proof resembles the strategy used in Example \ref{ex:ratsols1} and Example \ref{ex:ratsols2}: one considers the lines with rational slope passing through the given solution and proves that this gives infinitely many other solutions. This process does not help if we want to solve Diophantine equations in three variable and of degree greater than 2. Indeed we used the fact that, when intersecting the rational slope line with the degree two homogeneous curve, if one intersection point has rational coordinates then the other point of intersection has to be with rational coordinates too. If we move from degree 2 to, say, degree 3, than this assumption is not verified and, hence, new methods should be applied. For example in degree 3, one should consider more refined and involved techniques that make use of the Shimura-Taniyama-Wiles Conjecture which we are not going to discuss in the present article.

\subsection*{Pell's Equation}
In the previous examples we show ways to solve the problem of determining rational solutions to certain Diophantine equations of degree two. A related problem is to consider \emph{integral} solutions instead of rational ones. In this direction we are going to focus on one of the oldest Diophantine equations called the \emph{Pell Equation}. Such equation is usually written in the following form
\begin{equation}\label{eq:pell}
x^2-dy^2=1,
\end{equation}
with $d\in{\mathbb{N}}$ which we will assume is not a perfect square; the reason of this assumption is that if d were a perfect square than the right hand side would not be irreducible over $\Z$ and hence would have only trivial solutions (for a more extensive treatment of this equation we refer to \cite{LenstraPell}). We start by observing that if $(x,y)$ is a solution then $(\pm x, \pm y)$ is also a solution. This means that we can restrict ourselves to look for positive solutions to (\ref{eq:pell}).
We first address the problem of finding a particular solution to the Pell Equation. Note that a trivial solution always exists, namely $(1,0)$, but, for reasons we are going to discuss in a moment, we want to find another, non trivial, solution. This can be done by looking at the continued fraction expansion of $\sqrt{d}$.

\begin{example}\label{ex:sqrt5}
Consider the case with $d=5$, namely the Pell Equation
\begin{equation}\label{eq:pell5}
x^2 - 5 y^2 = 1.
\end{equation}
The continued fraction expansion of $\sqrt{5}$ is
\[
2 + \dfrac{1}{4 + \dfrac{1}{4 + \dfrac{1}{4 + \dots}}},
\]
and the convergents for this expansions are the following:
\[
\dfrac{2}{1},\ \dfrac{9}{4},\ \dfrac{38}{17},\ \dfrac{161}{72},\ \dots
\]
If we plug each fraction into (\ref{eq:pell5}), with $x$ being the numerator and $y$ the denominator we get
\begin{align*}
2^2 - 5 \cdot 1^2 &= -1 \\
9^2 - 5 \cdot 4^2 &= 81 - 80 = 1 \\
38^2 - 5 \cdot 17^2 &= 1444 - 1445 = -1 \\
161^2 - 5 \cdot 72^2 &= 25921 - 25920 = 1.
\end{align*}
We have just found \emph{two} non trivial solution to the Pell's Equation (\ref{eq:pell5})!
\end{example}

One can prove that the magic trick we just performed, i.e. looking at the continued fraction expansion of $\sqrt{d}$ and derive solutions to \ref{eq:pell5}, is indeed a theorem: one \emph{always} finds a non-trivial solution to Pell's Equation \ref{eq:pell} looking at the continued fraction expansion of $\sqrt{d}$. The solution with $\lvert x \rvert$ minimum is called the \emph{fundamental solution}. The reason why we are interested in such a solution is that starting from it we can build infinitely many others in the following way. Let $(x_1, y_1)$ be the fundamental solution of \ref{eq:pell}. Now every couple $(x_k, y_k)$ defined by the equality
\[
x_k + \sqrt{d} y_k = (x_1 + \sqrt{d} y_1)^k,
\]
is a solution of \ref{eq:pell} for all positive $k$. Therefore, the Pell's Equation has infinitely many solutions in positive integers
$x_k,y_k$. This fact is strictly related to Dirichlet's Theorem on units in quadratic number fields. If one considers the number field $\Q(\sqrt{d})$, then the left hand side of Pell's Equation can be factorized as follows:
\[
(x + y\sqrt{d})(x - y\sqrt{d}) = 1.
\]
This implies that $(x + y\sqrt{d})$ is a unit of norm 1. In this case Dirichlet's Theorem asserts that all these units form a finitely generated group of rank 1. The generator of this group is, in fact, related to the fundamental solution defined above.

\subsection*{Other Diophantine Problems}
We end this section with a couple of Diophantine problems which has played a major role in this subject, although they do not belong properly to classical Diophantine Equations. The first of this famous problems, which has been unsolved until 2002, is Catalan’s Conjecture which dates back to 1844. In its classical formulation it reads as follows:
\begin{conjecture}[Catalan’s Conjecture]
The equation
\begin{equation}\label{eq:catalan}
{x^p}-{y^q} = 1,
\end{equation}
where the unknowns $x, y, p$ and $q$ take integer values all $\geq2$, has only one solution,namely $(x, y, p, q) = (3, 2, 2, 3)$.
\end{conjecture}
Tijdeman’s 1976 result \cite{Tijdeman} showed that there are only finitely many solutions to Catalan's equation: more precisely, for any solution $x, y, p, q$ the number $\max\{p, q\}$ can be bounded by an effectively computable absolute constant. Once $\max\{p, q\}$ is bounded, only finitely many exponential Diophantine equations remain to be considered, and there are algorithms which complete the solution (based on Baker’s method). Such a bound has been computed, but it is somewhat large: Mignotte proved that any solution $x, y, p, q$ to (\ref{eq:catalan}) should satisfy $max\{p, q\} < 8 \cdot 10 16$ (see \cite{Mignotte} for a detailed historical discussion on this problem). Catalan’s claim was finally substantiated by P. Mih\u{a}ilescu. Note that the fact that the right hand side in Catalan’s equation is 1 is crucial: not so much is known if one replaces it by another positive integer.

Another very famous Diophantine problems is the so-called \emph{Waring's problem}. We present here only one of its many forms.

\begin{theorem}
Given an integer $k\geq2$, find the smallest integer $g(k)$ such that any non-negative integer can be represented
as a sum of $g(k)$ non-negative $k$-th powers.
\end{theorem}
Lagrange proved it for $g(2)=4$ and Hilbert showed that $g(k)$ is finite for all $k$. One of the major contributions towards the resolution of this problem was given by the Circle Method invented by Ramanujan, Hardy and Littlewood who improved the results in a critical way. The most difficult case turned out to be the one with $k=4$ which was finally solved by Balasubramanian, Deshouillers and Dress in \cite{Bala}, where they showed that $g(4)=19$.

\section{Diophantine Geometry - Rational and Integral points on curves}\label{sec:curves}
In the previous section we saw some problems related to Diophantine Equations for which we are interested in finding solutions, or, at least, in describing qualitatively their behavior. In fairly recent times a ``geometric'' approach has been introduced to study these problems and this theory goes under the name of \emph{Diophantine Geometry}.  One can describe Diophantine Geometry as the study of Diophantine Equations using geometric tools. In its simpler instance it asks for a ``geometric'' solution to the problem of finding reasonably simple necessary conditions for the solvability of
\[
f(x_1,\dots,x_n) = 0,
\]
where $f$ is a polynomial with coefficients in $\Q$, and $(x_1,\dots,x_n)$ is in either $\Z^n$ or $\Q^n$.
The term \emph{geometric} stands for the dependence on the properties of the complex algebraic variety defined by $f$, i.e. the set of \emph{complex} solutions to the previous equation.

More precisely and more generally, fixed a field $\kappa$ (which usually is either a number field, or a function field of an algebraic variety), we can describe Diophantine Geometry as the study of the set $X(\kappa)$ for an algebraic variety $X$ over $\kappa$. Rephrasing Lang, one can say that the main goal of Diophantine Geometry is to determine geometric properties of the quasi-projective variety $X: f = 0$ that characterize the set of points $X(\kappa)$, e.g. that imply that $X(\kappa)$ is either non-empty or finite or dense with respect to the Zariski topology. Here \emph{geometric} refers to properties that can be checked on the algebraic closure of $\kappa$.

\begin{example}\label{ex:dioph_ex}
Consider the equation
\[
\calC: -2x^2+3y^2-z^2=0.
\]
Since it is homogeneous, it defines naturally a projective conic $\calC_\C$ inside the projective plane $\PP^2_\C$, i.e. the set of complex solutions up to scalar multiplication. Now, as a projective variety over $\Q$, having the equation at least one solution over $\Q$, e.g. $\left[ 1: 1: 1\right]$, the conic $\calC$ is isomorphic to $\PP^1_\Q$, i.e. to the set
\[
\{ \left[ a: b\right] : a, b \in \Q \}.
\]
\begin{figure}
\begin{center}
\includegraphics[width=.5\textwidth]{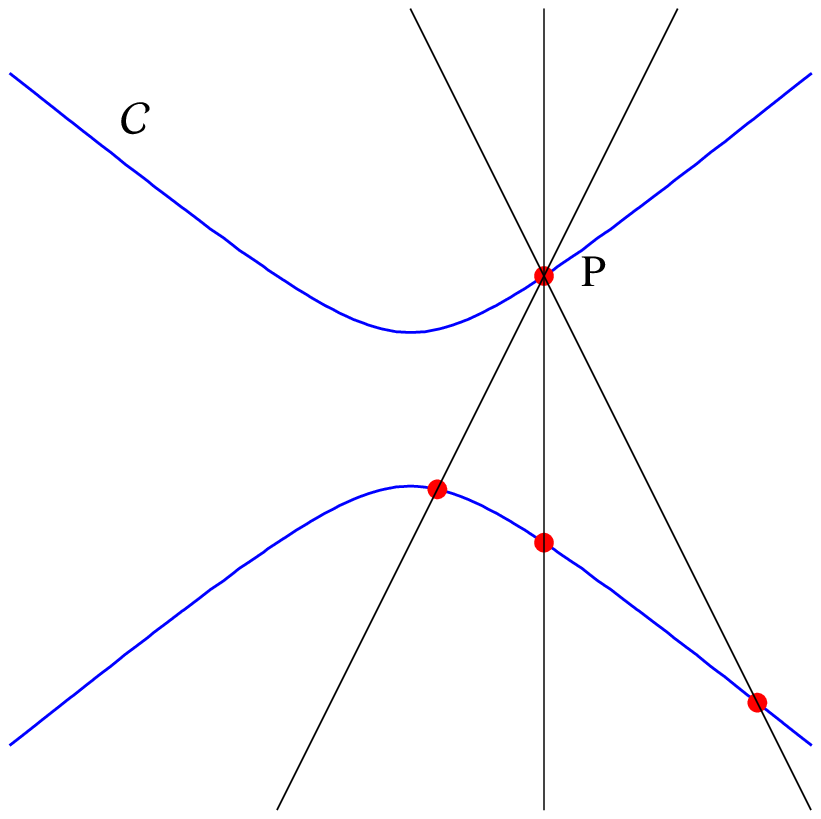}
\end{center}
\caption{A geometric description of an isomorphism between a conic with a rational point and $\PP^1$.}
\label{fig:P1}
\end{figure}
One can think of the isomorphism in the following way: consider $\PP^1$ as the set of lines in $\PP^2$ passing through a fixed point with rational coordinates, e.g. $P = \left[ 1:1:1 \right]$. Then, the map $\PP^1_\C \to \calC$ is defined by mapping each line $M$ through $P$ to the point $M \cap \calC \setminus P$. In the unique case in which $M$ is tangent to $\calC$ the image will be the point $P$ (see Figure \ref{fig:P1}). This set is in bijection with the solutions of the equation defining $L$ (for an explicit example of this procedure we refer to \cite{SilvermanTate}). Thus, since $\PP^1_\Q(\Q)$ is trivially an infinite set, there exist infinitely many solutions to the equation, i.e. infinitely many $\Q$-rational points on the line $L \subset \PP^2_\C$. Moreover, one can prove that this set is dense inside $L_\C \cong \PP^1_\C$ with respect to the Zariski topology. This gives some insight on how the knowledge of the \emph{geometry} associated to the Diophantine Equation gives information on the number and the behavior of the rational solutions without finding them explicitly as we did in Example \ref{ex:ratsols2}. 
\end{example}

In this section we will try to give some intuition for Diophantine Geometry problems in the one dimensional case, namely for affine and projective curves. By affine and projective curves we mean sets of solutions to equations of either the following types:
\[
f(x,y)=a, \qquad\qquad f(x_0,x_1,x_2) = 0,
\]
where $f$ is a polynomial which has integral (resp. rational) coefficients. Of course, affine equations $f(x,y)=a$ can be reduced to projective curves homogenizing the corresponding equation. We will assume that all the curves we are going to deal with are sufficiently \emph{nice}, i.e. every curve is geometrically irreducible and reduced; this is equivalent to require that the corresponding algebra
\[
\dfrac{\Q\left[ x_0, x_1, x_2\right]}{(f(x_0, x_1, x_2))}
\]
has neither zero divisors or nilpotents. Moreover we will assume that for every point of a projective curve there exists a unique tangent line, or, equivalently, that the curve is smooth over $\C$.

When considering complex solutions to a homogeneous equation in three variables, one is naturally led to consider the variety associated to it: it is a one dimensional complex projective variety that we assume to be smooth. What are the main geometric invariants
of such variety? From a birational point of view one could say that there exists only one discrete invariant called the \emph{genus}. The genus can be defined in several equivalent ways, e.g. it is the number of holes of the Riemann Surface associated to the algebraic curve or is the dimension of the $\Q$-vector space of regular differentials. For curves $\calC$ with ``nice'' singularities (ordinary) of degree $d$ the genus is given by the following useful formula:
\[
g(\calC) = \dfrac{(d-1)(d-2)}{2} - \sum_{P \text{sing}} \dfrac{m_P(m_P-1)}{2}.
\]
It turns out that for (smooth) projective curves defined over $\Q$ the genus discriminates the behavior of rational points.

\paragraph*{\textbf{genus $\mathbf{ = 0}$}}
From the genus formula given above, smooth curves of genus 0 are lines or conics. By Example \ref{ex:dioph_ex} we see that every time a genus 0 curve contains a rational point then it is isomorphic to $\PP^1$ over the rationals.
This implies that, for genus 0 curves $\calC$ with a point defined over $\Q$, the set $\calC(\Q)$ is infinite and dense in the Zariski Topology. However, there exist conics with \emph{no} rational points and which are \emph{not} isomorphic to $\PP^1$, e.g. $x^2+2y^2+z^2 = 0$. Summarizing this informal discussion one can state the following

\begin{proposition}\label{prop:g=0}
Given a smooth projective curve $\calC$ defined over $\Q$ of genus 0, one of the following holds:
\begin{enumerate}
\item If $\calC(\Q) \neq \emptyset$, then $\calC(\Q)$ is infinite and dense.
\item $\calC(\Q)$ is empty.
\end{enumerate}
\end{proposition}

In the second case, however, it is easy to see that allowing finite extensions $\kappa$ of the field $\Q$ one can always obtain a $\kappa$-rational point and, therefore, an isomorphism over $\kappa$ with $\PP^1$. In particular, even if the set of $\Q$-rational points can be empty, after a finite extension of the rationals one obtains infinitely many dense $\kappa$-rational points. For the conic $x^2+2y^2+z^2 = 0$, it is enough to consider the number field $\Q(i)$ with $i^2 = -1$. It is easy to see that by this extension the conic has $\Q(i)$-rational points, e.g. $\left[ i : 1 : i \right]$.

For integral points we look at affine curves defined by non-homogeneous equations and naturally embedded in an affine space $\Aff^n$. The genus can still be defined by looking at the completion of the curve, i.e. the projective curve associated to the homogenized equation. However, it could happen that this curve is not smooth, hence we should look at the so-called ``non-singular'' model associated to it. You can think of it as a procedure to obtain a smooth curve from a one which is singular (see \cite{Hartshorne} or \cite{Shafarevich} for explicit constructions). Then the genus of the affine curve can be defined to be the genus of the normalization of its completion. However, in the affine case, we have also another quantity that plays a role in the study of integral points, namely points at infinity. Let us give an explicit example:

\begin{example}\label{ex:integralg=0}
Consider the two affine curves in $\Aff^2$ given by the following equations:
\begin{align*}
\calC_1: f_1(x,y) &= x^2 - 5y^2 = 1, \\
\calC_2: f_2(x,y) &= x(1-x)y -1 = 0.
\end{align*}
Both curves have genus 0, the first having degree 2, and the second one having degree 3 and an ordinary point of multiplicity 2 (the point ``at infinity'' $\left[ 0: 1: 0 \right]$). By the previous discussion the corresponding completed curves, i.e. the set of zeros of the corresponding homogeneous equations, are isomorphic to $\PP^1$, since both have a rational point\footnote{We show this for smooth curves in Example \ref{ex:dioph_ex}; however the same construction can be carried over for $\calC_2$ using the lines passing through the double point.}. What can we say about the integral points? The equation defining the curve $\calC_1$ is a Pell equation; by Section 2, the set of integral solutions to this equation is infinite. Therefore, in our language, the set of integral points $\calC_1(\Z)$ is infinite too.
For the second curve however an easy computation shows that there are \emph{no} integral points! For seeing this it is sufficient to rewrite the equation defining $\calC_2$ as
\[
y = \dfrac{1}{x(1-x)}.
\]
Since $x$ should be an integer, in order for $y$ to be an integer both $x$ and $1-x$ should be $\pm 1$, which is impossible. We note passim that this is also related to the theory of $S$-unit equations, since finding $u=x$ and $v=1-x$ both invertible in $\Z$ is equivalent to solve the equation $u+v=1$ in units for $\Z$.

\noindent How can we explain this completely opposite behavior of the sets $\calC_1(\Z)$ and $\calC_2(\Z)$ while both curves have completion isomorphic to $\PP^1$? The answer is that we have to take into account how ``far'' this affine curves are from $\PP^1$, i.e. how many points should we add to $\calC_1$ and $\calC_2$ to get $\PP^1$.
For $\calC_1$ one sees immediately that the point at infinity, given by the solutions of
\[
\begin{cases}
x^2 - 5y^2 -z^2 = 0 \\
z = 0
\end{cases}
\]
are $\{ \left[\sqrt{5}:1:0\right], \left[-\sqrt{5}:1:0\right]\}$. For the curve $\calC_2$ one can prove that such curve is isomorphic (possibly after enlarging the base field) to $\PP^1$ minus three points, e.g. $1,0,\infty$. This can be seen for example by comparing the ring of regular functions on the two affine curves. In particular $\calC_2$ has three points at infinity compared to the two of $\calC_1$. This is a crucial feature that will explain the different behavior of the sets $\calC_1(\Z)$ and $\calC_2(\Z)$.

\end{example}

From the previous example one can see how for genus zero curves one should take into account also the number of points at infinity. In this direction, the distribution of integral points on affine curves of genus zero is governed by this result which is a corollary of a deep theorem of Siegel:

\begin{theorem}[Siegel]\label{th:Sigelg0}
Let $\calC$ be an affine irreducible curve defined over $\Q$. Let $\tcal{C}$ be the completion of its normalization and $g = g(\tcal{C})$ its (geometric) genus. Suppose that $g=0$ and $\tcal{C}\setminus \calC$ contains at least three points. Then $\calC(\Z)$ is a finite set in $\Aff^m(\Z)$. If the set $\tcal{C}\setminus \calC$ contains less than three points, at most after a finite extension of $\Z$ the set $\calC(\Z)$ is infinite.
\end{theorem} 

\paragraph*{\textbf{genus $\mathbf{ = 1}$}}
Smooth projective curves of genus one with a specified base point are \emph{elliptic curves} and are among the most studied algebro-geometric objects, having a natural group structure that enriches their arithmetic and their geometry. For the problem we are concerned in this survey (namely a description of the set of integral and rational points) the case of elliptic curves is reduced to two very deep theorems, namely Siegel's Theorem and Mordell's Theorem, respectively. We have just seen Siegel's Theorem in the case of genus zero for integral points on rational curves and we will see at the end of this section its general formulation, that encompass all the possible genus and number of points at infinity. We will begin, as before, from the problem of describing the set of rational points on an elliptic curve.

Let $E$ be a smooth projective curve of genus 1: the geometric picture of this curve in the projective plane over the complex numbers is a torus. This follows from the fact that one can construct an elliptic curve by looking at the quotient of $\C$ by a lattice $\Lambda$ and embed this quotient via a map defined using the Weierstrass $P$-function. Such a map can be proved to be an isomorphism of groups: this fact in particular implies that an elliptic curve constructed in this way is naturally an abelian group. It follows also that the group structure gives a natural choice for the base point, namely the identity element of the group.

\begin{figure}
\begin{center}
\includegraphics[width=.5\textwidth]{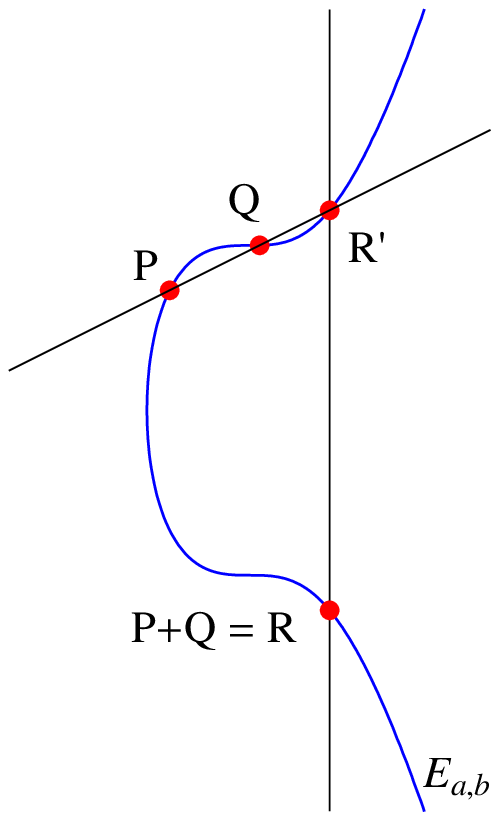}
\end{center}
\caption{The law of addition in an elliptic curve.}
\label{fig:Elliptic}
\end{figure}

From a more algebraic point of view, elliptic curves over $\Q$ are curves defined by equations of the type
\begin{equation}\label{eq:ell}
E_{a,b}: y^2 = x^3 + a x + b,
\end{equation}
where $a$ and $b$ are rational numbers such that the discriminant
\[
\Delta_{a,b} = -16(4a^3 + 27b^2)
\]
is non-zero. This is equivalent to require that the algebraic curve is smooth (as we saw in Example \ref{ex:integralg=0}, if the curve is defined by an equation of degree 3 but it possesses a singular double point, then its genus is 0; this implies that the curve is not elliptic). One also requires that there exists at least one rational solution of (\ref{eq:ell}) usually this point is the point at infinity, i.e. the unique solution to (\ref{eq:ell}) obtained by homogenizing $E_{a,b}$ using the extra variable $z$ and then imposing $z=0$; in this case the only (flex) point at infinity is the point $\left[ 0: 1: 0\right]$. As already mentioned above the set of rational points on such curves possesses a natural group structure whose operation can be seen geometrically as it follows. Let $I$ be the flex (rational) point at infinity, whose existence is required in the definition of elliptic curve. Given any two rational points $P,Q \in E_{a,b}$, we define their sum $R := P + Q$ in the following way: first let the point $R'$ be the third point of intersection of the line passing through $P$ and $Q$ with $E_{a,b}$; then $R$ is the point of intersection of the line joining $R'$ and $I$ with the cubic. A picture of this procedure can be seen in Figure \ref{fig:Elliptic}.

One can extend this definition also to the case in which $R'$ coincides with $P$ or $Q$ and prove that the point $I$ acts as the identity for this operation. Moreover, the addition just defined can be proven to be an associative and commutative operation on the set of rational points $E_{a.b}(\Q)$. This implies that $(E_{a,b}(\Q),+)$ is an abelian group. What kind of group? To answer this question Mordell proved his celebrated theorem which reads as follows:

\begin{theorem}[Mordell \cite{Mordell1922}]\label{th:MordellEl}
Given a non-singular elliptic curve $E_{a,b}$, the set of rational points $E_{a,b}(\Q)$ is a finitely generated abelian group.
\end{theorem}

In particular, the theorem implies that at most after a finite extension of $\Q$, the set or rational points on an elliptic curve is an infinite set. For a detailed description of the proof of Theorem \ref{th:MordellEl} we refer to \cite{SilvermanTate}.
What can be said for integral points? As for the genus 0 case the result below following from Siegel's Theorem gives a definitive answer.

\begin{theorem}[Siegel]\label{th:Sigelg1}
Let $\calC$ be an affine irreducible curve defined over $\Q$. Let $\tcal{C}$ be the completion of its normalization and $g = g(\tcal{C})$ its (geometric) genus. Suppose that $g=1$ and $\tcal{C}\setminus \calC$ contains at least one point. Then $\calC(\Z)$ is a finite set in $\Aff^m(\Z)$.
\end{theorem} 

Again we see how the genus and the number of points at infinity completely describe the distribution of integral points. In particular one sees that, after a finite extension of $\Q$ there are always infinitely many integral points on an elliptic curve. This follows from the fact that the equation defining points at infinity has always a solution at most up to a finite extension of the base field. This implies that, up to this extension, there exists always a point at infinity for an elliptic curve which implies that, by Theorem \ref{th:Sigelg1}, the set of integral point is an infinite set.

\paragraph*{\textbf{genus $\mathbf{ \geq 2}$}}
The study of rational points on higher genus curves has been historically one of the most challenging mathematical problems and it is related to the well known Mordell Conjecture. The history of the Mordell Conjecture starts with the article \cite{Mordell1922} by Mordell. In this seminal paper, after proving Theorem \ref{th:MordellEl}, he states one of the most famous problems in Arithmetic Geometry, which was proved by Faltings only sixty years after the original formulation. The Conjecture, now Faltings' Theorem, reads as follows.
\begin{theorem}[Faltings \cite{Faltings1983}, Mordell Conjecture]\label{th:MordellFaltings}
Let $K$ be a number field and let $\calC$ be a curve defined over $K$ of genus greater than one. Then $\calC(K)$ is a finite set.
\end{theorem}
This statement was highly non-trivial and only some particular case were known at Mordell's time. Many mathematicians, although recognizing its power, were not convinced by the conjectured result. Andr\'e Weil commented
\begin{verse}
Nous sommes moins avanc\'{e}s \`{a} l'\'{e}gard de la Conjecture de Mordell. Il s'agit l\`{a} d'une question qu'un arithmeticien ne peut gu\`{e}re manquer de se poser; on n'aper\c{c}oit d'ailleurs aucun motif s\'{e}rieux de parier pour ou contre. \cite{Weil1979} \vspace{.3cm}\\
We are less advanced in respect of the Mordell Conjecture. This is a problem that every arithmetician can hardly not ask himself; nevertheless we do see no serious reason to bet for or against its truth.
\end{verse}
In his general audience exposition of Mordell Conjecture and Falting's ideas \cite{Bloch1984}, Spencer Bloch wrote \textit{Probably most mathematicians would have agreed with Weil (certainly I would have) until [...] a German mathematician, Gerd Faltings, proved the Mordell Conjecture}. This is even more revealing, if we take into account the proof by Grauert \cite{Grauert1965} and Manin \cite{Manin1963} (although with a gap pointed out and corrected by Coleman \cite{Coleman1990}) of the function field case.

Nevertheless, in 1983 Faltings presented a proof of Theorem \ref{th:MordellFaltings} as a consequence of his proof of Tate Conjecture and Shafarevich Conjecture. His argument uses very refined and difficult tools, like Arakelov Theory on moduli spaces, semistable abelian schemes and $p$-divisible groups. Vojta in \cite{Vojta1991} (and previously for function fields in \cite{Vojta1989}) gave another proof which uses ideas from classical Diophantine approximations together with technical tools of intersection theory on arithmetic threefolds developed by Gillet and Soul\`e. After this new proof, Faltings in \cite{Faltings1991} gave another simplification, eliminating the use of Riemann-Roch Theorem for arithmetic threefolds: using his new ideas he was able to extend previous results and to prove a conjecture formulated by Lang. Another simplification of both Vojta and Faltings' proofs was given by Bombieri in \cite{Bombieri1990}, combining idea from Mumford \cite{mumford1965} together with the ones in the aforementioned papers.

A full proof of Theorem \ref{th:MordellFaltings} goes beyond the scope of this survey, therefore we refer to the following books that contains a detailed and comprehensive discussion of the original proofs, together with their subsequent simplifications: Bombieri and Gubler \cite{BombieriGubler} and Hindry and Silverman \cite{HindrySilverman} discuss Bombieri's approach to Theorem \ref{th:MordellFaltings}. For an exposition of the ideas of Faltings' original paper the main source is Faltings and W{\"u}stholz notes \cite{Faltings1986}; another exhaustive treatment of the original proof together with its link to Tate Conjecture and Shafarevich Conjecture can be found in Zarhin and Parshin article \cite{Zarhin2009}\footnote{For bibliographic history of this article we refer to the summary in the article's arXiv page. We just notice that the paper first appeared as an appendix to the Russian version of Lang's book \emph{Fundamentals of Diophantine Geometry}.}.

As for integral points, as you may expect, the Siegel's Theorem gives the full answer. We state here the original version of Siegel's Theorem encompassing all possible genus and number of points at infinity, which reads as follows:

\begin{theorem}[Siegel 1929,\cite{Siegel1929}]\label{th:SiegelNF}
Let $\calC$ be an affine irreducible curve defined over $\Q$ and embedded in an affine space $\Aff^m$. Let $\tcal{C}$ be the completion of its normalization and $g = g(\tcal{C})$ its (geometric) genus. Suppose that either $g \geq 1$ or $g=0$ and $\tcal{C}\setminus \calC$ contains at least three points. Then $\calC(\Z)$ is a finite set in $\Aff^m(\Z)$.
\end{theorem}

\begin{remark}\label{rmk:SiegelTh1}
Historically, Siegel's original presentation of the result was split into two parts: one dealing with the case of genus zero and three points at infinity and the remaining one dealing with genus greater or equal than one.
The theorem has been extended to number fields by Mahler in \cite{Mahler1934} for genus 1. The result has been finally extended to arbitrary finite set of places by Lang in \cite{Lang1960} using an extension of Thue-Siegel-Roth Theorem by Ridout in \cite{Ridout1958}.
\end{remark}

The original proof of Siegel's Theorem uses the so-called Thue-Siegel Theorem from Diophantine Approximation together with some properties of theta characteristics. The modern version of this proof relies on the more general Roth Theorem and on the theory of heights in the Jacobians. In \cite{Corvaja2002} Corvaja and Zannier gave another proof of Theorem \ref{th:SiegelNF} avoiding the embedding in the jacobian and replacing the use of Roth's Theorem with the stronger Schmidt Subspace Theorem. The importance of this new reformulation, left aside the fact that it simplifies Siegel's argument, relies on an extension to higher dimensions which will be important in the next section.

The main importance of Faltings' Theorem (together with Siegel's Theorem \ref{th:SiegelNF}) for the purpose of this article is the following corollary which completely describes the distribution of integral and rational points on algebraic curves:

\begin{theorem}[Arithmetic classification of curves]\label{th:arith_class_curves}
Let $\calC$ be a projective, geometrically irreducible and non-singular curve defined over $\Q$. Then, at most after a finite extension of $\Q$\footnote{This assumption is made in order to give a unified treatment of the case of genus zero and one. Indeed, at most after a finite extension, every algebraic curve possess a rational point.}, the following description holds:
\begin{center}
\begin{tabular}{|l|c|l|c|}
\hline
Genus		&Rational points 		& Points at infinity 		& Integral points \\
\hline
$g=0$		& Infinite set				& \qquad $\leq 2$								& infinite set \\
\hline
$g=0$		& Infinite set				& \qquad $\geq 3$								& finite set \\
\hline
$g=1$		& Fin. generated group				& \qquad $=0$				& infinite set \\
\hline
$g=1$		& Fin. generated group				& \qquad $\geq 1$		& finite set \\
\hline
$g \geq 2$& Finite set				& \qquad Arbitrary								& finite set \\
\hline
\end{tabular}

\end{center}
\end{theorem}

The previous table gives the final answer to the question whether there exist infinitely many rational or integral points in an algebraic curve. This can be translated into a strategy useful to deal with arbitrary polynomial equations in two or three variables with integral or rational coefficients: namely given an equation $f=0$ one considers the normalization of the homogeneous equation $f_{\text{hom}}=0$ associated to it and computes the genus of such non-singular projective curve (at most after normalizing it) and the number of points at infinity; then Theorem \ref{th:arith_class_curves} tells whether the sets $f=0$ over $\Z$ and $f_{\text{hom}}=0$ over $\Q$ are finite or not.

\section{Rational and Integral points on surfaces}\label{sec:surf}
In this last section we are going to discuss the description of rational and integral points on algebraic varieties of (complex) dimension 2, or equivalently the description of solution to polynomial equations in four (projective) variables. As already mentioned in the introduction, we are going to need more tools from Algebraic Geometry. Therefore this section will require more prerequisites than the others. Moreover, since this does not require too much work, we are going to move from $\Q$ to number fields, i.e. allowing finite extensions of the rational numbers, and considering the correspondent notion of integers, namely $S$-integers and units with respect to a finite set of places $S$.

From Theorem \ref{th:MordellFaltings} one can see how the geometric properties encoded by the genus \emph{govern} the arithmetic of the curve. Seeking a generalization to higher dimensions, and in particular to surfaces, one is led to study which geometric features of the underlying complex variety determine the distribution of rational points. At the same time one should carefully consider whether the same questions arising in the one-dimensional case can be extended to surfaces. Indeed one could start by studying whether rational points on higher dimensional varieties are finite or not. However, as the following example shows, being finite could be a tricky property to check and does not fully describe the distribution of rational and integral points on surfaces. 

\begin{example}[Corvaja and Zannier, Turchet]\label{ex:cubic}
Let $\tilde{X}$ be a smooth cubic surface defined over $\Q$ and let $H_1, H_2$ be two hyperplane sections defined over $\Q$ such that $H_1\cup H_2$ consists of 6 lines. Corvaja and Zannier in \cite{Corvaja2009a} proved that the set of integral points on $X = \tilde{X} \setminus (H_1 + H_2)$ is not Zariski dense. On the other hand, one can prove (see \cite{TurchetMaster}) that the only families containing infinite integral points are the 21 remaining lines in $\tilde{X}$ (it is a well known and beautiful fact that each cubic surface contains exactly 27 lines and for the discussion of Section \ref{sec:curves} we know that each line possesses infinitely many rational points, at most after a finite extension of $\Q$).

This example shows how in a complement of two completely reducible hyperplane sections in a smooth cubic the integral points are "almost finite" in the sense that, removing a finite number of subvarieties (or a proper subvariety consisting of the union of these), the integral points are finite in the surface. In particular, the closure of the set of integral points is a proper subvariety of the affine surface. Therefore from an informal point of view we would like to consider such a variety as one of those who possesses ``few'' integral points.
\end{example}

Following the idea of Example \ref{ex:cubic}, we will focus on the problem of determining whether the set of rational or integral points in an algebraic surface is dense or not with respect to the Zariski topology. In order to extend Mordell and Faltings ideas to surfaces, we have to look for geometric properties of algebraic surfaces which could imply that the set of rational points is not dense on the surface. Therefore, we need geometric properties replacing or, better, extending the role played by the genus in dimension 1. With this goal in mind we recall the following

\begin{definition}[Kodaira dimension]\label{def:KodDim}
Let $X$ be a smooth projective algebraic variety and let $K_X$ be one if its canonical divisors. For each $m\geq 1$ such that the pluricanonical linear system $\lvert m K_X \rvert$  is not empty, i.e. such that $h^0 (X, \calO(m K_X)) \neq 0$, let
\[
\Phi_{m K_X} : X \to \PP^N
\]
be the associated map. The \emph{Kodaira dimension} of $X$ is defined to be the number
\[
\kappa(X) = \begin{cases} 
							-1 \quad & \text{ if } h^0 (X, \calO(m K_X)) = 0\quad \forall m, \\
							\max \dim \Phi_{m K_X} (X) & \text{ otherwise }.
					\end{cases}
\] 
\end{definition}

Note that the Kodaira dimension for curves is given by $\kappa = \min \{1, g - 1\}$. In particular, we could rephrase the description given in Section \ref{sec:curves} in the following way

\begin{corollary}\label{cor:KodCurves}
Given a smooth projective curve $\calC$ defined over a number field $K$, the set of rational points $\calC(K)$ is not dense if and only if $\kappa(\calC) = 1 = \dim \calC$.
\end{corollary}

Motivated also by the previous result, we recall the following terminology that extends the property of having genus greater than 1 for curves.

\begin{definition}[General type varieties]\label{def:GenType}
Let $X$ be an algebraic variety. If $\kappa(X) = \dim X$, then $X$ is said to be of \emph{general type}.
\end{definition}

The idea behind Lang and Bombieri conjecture on distribution of rational points on surfaces is that a rough analogous of the behavior exhibited by algebraic curves could hold also for surfaces. First of all, we recall the Kodaira classification of surfaces, which reads as follows:

\begin{theorem}[Kodaira Classification of Surfaces]\label{th:KodClassification}
Let $X$ be an algebraic smooth surface and let $\kappa = \kappa(X)$ its Kodaira dimension. Then the following classification holds:
\begin{itemize}
\item $\kappa = -1$: $X$ is either a Rational or a Ruled surface.
\item $\kappa = 0$: $X$ belongs to one of the following four classes: Abelian, hyperelliptic (or bi-elliptic), K3 or Enriques.
\item $\kappa = 1$: $X$ is an Elliptic Surface.
\item $\kappa = 2$: by definition $X$ is of General Type.
\end{itemize}
\end{theorem}

We are interested in the behavior of the set of rational points for each family of surfaces listed in Theorem \ref{th:KodClassification}. However, differently from the curves' case, for surface most of the results are still conjectural. Let us look more closer to each item appearing in the list:

\begin{itemize}
\item Let us consider the first case: each Rational or Ruled surface defined over a number field $K$ is covered by rational curves which, by Theorem \ref{th:MordellFaltings} have infinitely many rational points. Therefore for all the surfaces in this class the set of $K$-rational points is (potentially) dense.
\item The case of null Kodaira dimension is more involved: it is known that rational points are potentially dense, i.e. dense at most after a finite extension of the ground field, for abelian varieties and for Enriques surfaces \cite{BogomolovTschinkel98}. There are several proved results of density of rational points for some classes of K3 surfaces \cite{BogomolovTschinkel00} and for Hyperelliptic surfaces \cite{BogomolovTschinkel99}. Conjectures predicts that for each of these classes the rational points \emph{are} potentially Zariski dense.
\item For elliptic surfaces of Kodaira dimension one there is a gap for an arithmetic classification of surfaces based solely on Kodaira dimension. In fact one can easily construct example of surfaces with $\kappa = 1$ that posses either a potentially dense set of rational points or a non-dense one. Consider two fibrations $X \to \calC$ defined over a number field $K$ having elliptic curves as fibers; suppose that the genus of the base curve $\calC$ is greater than one: then from Faltings' Theorem \ref{th:MordellFaltings} $X$ has a non (potentially) dense set of rational points. On the other hand if $\calC = \PP^1$ and there exists infinitely many sections over (a finite extension of) $K$, $K$-rational points are potentially dense in $X$. In both cases it may happen that $X$ has Kodaira dimension equal to one.
\item For surfaces of general type it is expected that the set of rational points is not potentially dense: this has been conjectured independently by Bombieri and Lang. Bombieri addressed the problem of degeneracy of rational points in varieties of general type in a lecture at the University of Chicago in 1980, while Lang gave more general conjectures centered on the relationship between the distribution of rational points with hyperbolicity and Diophantine approximation (see \cite{Lang1997} and \cite{Lang1974}). The conjecture reads as follows:
\end{itemize}

\begin{conjecture}[(Weak) Bombieri-Lang]\label{conj:BL}
Let $X$ be a surface of general type defined over a number field $K$. Then the set of $K$-rational points of $X$ is not Zariski dense.
\end{conjecture}

Evidences for Bombieri-Lang Conjecture come from the following conjecture due to Lang and proved by Faltings in \cite{Faltings1991} and \cite{Faltings1994}.

\begin{theorem}[Lang Conjecture - Faltings' Big Theorem]\label{th:FaltingsLang}
Let $A$ be an abelian variety over a number field $K$ and let $X$ be a geometrical irreducible closed subvariety of $A$ which is not a translate of an abelian subvariety over $\Kbar$. Then $X \cap A(K)$ is not Zariski dense in $X$.
\end{theorem}

See \cite{Hindry98} for a detailed introduction and explanation of this conjecture. From the previous theorem it follows a corollary which gives several evidences to Bombieri-Lang;

\begin{corollary}\label{cor:gen_type_in_abelian}
If $X$ is a smooth projective variety of general type defined over a number field contained in an abelian variety, then the set of rational points of $X$ is not Zariski dense.
\end{corollary}

Following Noguchi's proof \cite{Noguchi1981} in the function field case for varieties whose cotangent bundle is ample (which implies that the variety is of general type) and using Faltings' Big Theorem (cfr. Theorem \ref{th:FaltingsLang}), Moriwaki in \cite{Moriwaki1995} obtained another evidence for Bombieri-Lang Conjecture which reads as follows:

\begin{theorem}[Moriwaki]\label{th:Moriwaki}
Let $X$ be a projective variety over a number field $K$. If the sheaf of differentials $\Omega^1_{X/K}$ of $X$ over $K$ is ample and generated by global sections, then the set of $K$-rational points of $X$ is finite.
\end{theorem}

Other evidences for Bombieri-Lang Conjecture comes from related examples and conjectures for the distribution of rational curves in general type surfaces, such as Bogomolov Theorem on the finiteness of rational and elliptic curves on general type surfaces with $c_1^2 > c_2$ \cite{Bogomolov1978}.

In the same way Conjecture \ref{conj:BL} extends Faltings' Theorem \ref{th:MordellFaltings}, it is natural to ask whether a similar extension exists for Siegel's Theorem \ref{th:SiegelNF}. The answer is positive and it is related to Vojta's ''landmark Ph.D. Thesis", which gave the basis for a systematic treatment of analogies between Nevanlinna Theory and Diophantine Geometry over number fields. Based on this analogy Vojta formulated a set of far-reaching conjectures. For a detailed description we refer to Vojta's papers \cite{Vojta1987} and \cite{Vojta2011} as well as chapters in the books \cite{HindrySilverman} and \cite{BombieriGubler}.

In order to properly stating the general Vojta's Conjecture we would have needed all the theory of Global and Local Heights on algebraic varieties which goes beyond the aims of this survey. Instead, we are going to state a weaker version, which uses ideas of Lang, where the condition of being of general type in Conjecture \ref{conj:BL}, is replaced by the condition of being of \emph{log-general type}. We need the following

\begin{definition}\label{def:logGT}
Let $X$ be a smooth projective variety and $D$ a normal crossing divisor on $X$. $X$ is said to be of \emph{logarithmic general type}, or log-general type, if $K_X + D$ is big for a canonical divisor $K_X$ of $X$.
\end{definition}

Using this definition we can state the following

\begin{conjecture}[Lang-Vojta]\label{conj:LVNNii}
Let $X$ be a quasi projective surface of \emph{log-general type} defined over a number field $K$ and let $\Oint{S}$ be the ring of $S$-integers for a finite set of places of $K$ containing the archimedean ones. Then the set $X(\Oint{S})$ is not Zariski dense.
\end{conjecture}

From the statement, one could immediately see how, in dimension one, the Conjecture encompasses Theorem \ref{th:SiegelNF}.
One important case of Conjecture \ref{conj:LVNNii} follows directly from Schmidt's Subspace Theorem:

\begin{corollary}\label{cor:4linesNF}
Given 4 lines $D_1,\dots,D_4$ in $\PP^2_K$ defined over a number field $K$ and $S$ a finite set of places containing the archimedean ones, the set of $S$-integral points on the complement of $D_1 + \dots + D_4$ is not Zariski dense.
\end{corollary}

Notice that the divisor formed by four lines in general position, makes the complement $\PP^2 \setminus D$ where $D = D_1 + \dots D_4$ a variety of log-general type, because $K_{\PP^2} + D \sim \calO_{\PP^2}(1)$ is an ample divisor. In particular Lang-Vojta Conjecture \ref{conj:LVNNii} holds for the complement of at least $4$ hyperplanes in general position in $\PP^2$.

In more recent times, using their new proof of Siegel's Theorem, Corvaja and Zannier obtained a number of strong results on degeneracy of integral points on surfaces by means of this new strategy. Among these results we cite the following which extends Corollary \ref{cor:4linesNF} to more general situations.

\begin{theorem}[\cite{Corvaja2006}]\label{CZNF}
Let $X$ be a geometrically irreducible nonsingular projective surface defined over a number field $K$ and let $D_1,\dots,D_4$ be irreducible effective divisors such that
\begin{enumerate}
\item No three of them shares a common point;
\item For all $i \neq j$, $\supp D_i \cap \supp D_j \neq \emptyset$;
\item For all $i \neq j$, $D_i \underset{\text{num}}{\sim} m D_j$ for a certain $m = m_{i,j} \in \Z$.
\end{enumerate}
Then no set of $S$-integral points in $X \setminus D$ is Zariski dense.
\end{theorem}

\begin{remark}
\begin{itemize}
\item Clearly for $X = \PP^2$ all the hypotheses of the Theorem are verified for the divisor $D$ consisting of four lines in general position and hence the Theorem implies Corollary \ref{cor:4linesNF}. In particular each $D_i=L_i$ is effective (even ample in this case), the general position hypothesis implies condition 1 and 2, while condition 3 follows from the fact that, being $\Pic \PP^2$ of rank 1, all classes of lines are linearly equivalent. 
\item The proof of the previous Theorem relies on a generalization of the ideas of \cite{Corvaja2002} where a suitable choice of a linear system of multiples of the irreducible components $D_i$ replaces the linear spaces of rational functions with prescribed order of zeros at limit points.
\end{itemize}
\end{remark}

A more general theorem can be found in \cite{Corvaja2004} and, with some modification and extension, in \cite{Corvaja2006}. It is worth mentioning also a corollary obtained by Levin in \cite{Levin2009} where he was able to drop the third hypothesis on the $D_i$ provided that $D_i$ is ample for every $i$.

As one can see from the previous theorems, a lot of results are known when the divisor at infinity $D$ has many irreducible components. This comes form the fact that in these cases one can reduce the problem to a Diophantine Approximation problem and reduce the problem to one which can be solved with an application of the Subspace Theorem. On the other hand, for irreducible $D$ only few partial results are known.

We should mention that a function field version of Conjecture \ref{conj:LVNNii} has been extensively studied in the last years, and all the results mentioned above have been proved also in this case. Moreover, the techniques available for function field arithmetic, have allowed to obtain broader and deeper results. In particular, for $\PP^2$, Corvaja and Zannier proved in \cite{Corvaja2005} the function field version of Lang-Vojta Conjecture for complements of divisors $D$ of degree at least 4 and with at least three components. At the same time, moving from Hyperbolicity problems related to Kobayashi Conjecture, Xi Chen in \cite{Chen2001} and \cite{Chen2004}, and Pacienza Rousseau in \cite{Rousseau2007}, proved that the conjecture holds for complements of \emph{very general} divisors of degree at least 5, without any limitation on the number of irreducible components. Recently, in his Ph.D. Thesis \cite{TurchetPhD}, the second author extended this results proving that the conjecture holds for all affine surfaces $\PP^2\setminus D$ of log general type, provided that $D$ has simple normal crossing and it is very general.

\bibliographystyle{alphanum}	
\bibliography{Mendeley}

\end{document}